# WHEN CAN THE TWO-ARMED BANDIT ALGORITHM BE TRUSTED?

By Damien Lamberton, Gilles Pagès and Pierre Tarrès

Université de Marne-la-Vallée, Université Paris 6 and Université Paul Sabatier

We investigate the asymptotic behavior of one version of the so-called two-armed bandit algorithm. It is an example of stochastic approximation procedure whose associated ODE has both a repulsive and an attractive equilibrium, at which the procedure is noiseless. We show that if the gain parameter is constant or goes to 0 not too fast, the algorithm does fall in the noiseless repulsive equilibrium with positive probability, whereas it always converges to its natural attractive target when the gain parameter goes to zero at some appropriate rates depending on the parameters of the model. We also elucidate the behavior of the constant step algorithm when the step goes to 0. Finally, we highlight the connection between the algorithm and the Polya urn. An application to asset allocation is briefly described.

**Introduction.** The aim of this paper is to deeply investigate the asymptotic behavior of the so-called two-armed bandit algorithm. This stochastic approximation procedure is widely known in the fields of mathematical psychology and learning automata (see [13] and [15]). Our own motivations are both theoretical and practical as it will be seen further on. Let us first introduce the algorithm itself in a financial context, namely as an adaptive optimal asset allocation model.

Imagine a fund managed by only two traders, say $A$ and $B$: every day each of them is in charge of a percentage of the fund, which may vary from day to day. The few wealthy investors (the shareholders) who created the fund wish ideally to allocate the whole fund to the most efficient trader, but of course they do not know who he is. They simultaneously want to make some advantage of the performances of the best trader as soon as possible. This means they need to devise a periodic re-allocation procedure of the fund to the traders based on their (daily or monthly) performances. On the other hand, this procedure should be not too "upsetting" to the traders in order to preserve their motivation and self-confidence: one way is to enhance reward rather than punishment. Taking all these specifications into account suggests to proceed as follows: let $X_n$ be the fraction of the fund managed by trader $A$ during day $n$, the fraction $1 - X_n$ being managed by trader $B$. Every day, one trader is chosen *at random* and his performances of the day are evaluated. Assume it is $A$ for a moment. If they are considered as outstanding, trader $A$ is rewarded by an extra-allocation for day $n + 1$ of $\gamma_{n+1}$ *times* the fraction managed by trader $B$ during day $n$ (whatever the performances of trader $B$ are since he was not checked). So trader $A$ will manage a fraction $X_n + \gamma_{n+1}(1 - X_n)$ of the fund during the day $n + 1$. If his performances are not high enough to deserve a reward, nothing happens. The same procedure is applied to trader $B$ when he is checked: if $B$ has outstanding performances, he is awarded an extra allocation $\gamma_{n+1}$ *times* the share managed by $A$ during day $n$ so that, during day $n + 1$, the share managed by $A$ will be reduced to $X_{n+1} = X_n - \gamma_{n+1} X_n$ (whatever his performances on day $n$ were). One models the daily performance evaluations of $A$ and $B$ by two sequences of events $(A_n)_{n \geq 1}$ and $(B_n)_{n \geq 1}$, respectively: $A_n = \{A$'s performances on day $n$ are outstanding$\}$ and $B_n = \{B$'s performances on day $n$ are outstanding$\}$.

A natural policy for the investors of the fund is to reduce the risks induced by this strategy by controlling the largest possible part (in average) of the whole fund. So tossing up for the checked trader with a fair coin is not appropriate. What seems more efficient is to use for the daily toss a biased (virtual) coin so that the probability for trader $A$ or $B$ to be checked at the end of day $n$ is equal to the share of the fund they managed that day, namely $X_n$ and $1 - X_n$, respectively. This virtual coin can be tossed by generating on a computer some i.i.d. random numbers $U_n$, $n \geq 1$ and by setting

$$\{A \text{ is checked at the end of day } n\} = \{U_{n+1} \leq X_n\},$$
$$\{B \text{ is checked at the end of day } n\} = \{U_{n+1} > X_n\}.$$

All this leads to the following dynamics for $X_n$: for every $n \geq 0$,

$$(1) \quad X_{n+1} = X_n + \gamma_{n+1}((1 - X_n)\mathbf{1}_{\{U_{n+1} \leq X_n\} \cap A_{n+1}} - X_n \mathbf{1}_{\{U_{n+1} > X_n\} \cap B_{n+1}}),$$
$$X_0 = x \in [0, 1],$$

where $(\gamma_n)_{n \geq 1}$ is the sequence of gain parameters (or *steps*) satisfying

$$(2) \quad \forall n \in \mathbb{N}^*, \ \gamma_n \in (0, 1) \quad \text{and} \quad \Gamma_n := \gamma_1 + \cdots + \gamma_n \to +\infty \quad \text{as } n \to +\infty.$$

[Note this includes the constant step setting $\gamma_n = \gamma \in (0, 1)$.] The fact that $\gamma_n$ lies in $(0, 1)$ is induced by the modelling (it is a percentage). On the other hand, the fact that $\Gamma_n$ goes to infinity is a necessary condition to "forget" the starting value: if $\lim_n \Gamma_n < +\infty$, $X_n$ would still converge a.s. toward a random variable $X_\infty$, but one could not show that $X_\infty$ takes its values in $\{0,1\}$.

This recursive random procedure was first introduced by Norman in mathematical psychology (see [13]) and then, independently, by Shapiro and Narendra in the engineering literature as a linear learning automata (see [15]). In this field it is known as the Linear Reward-Inaction ($L_{R-I}$) scheme (see the survey [10] and the book [11] by Narendra and Thathachar about learning automata theory). In both cases, only the constant step setting is considered. The application to optimal adaptive asset allocation in a financial context has been developed in [12].

The algorithm (1) is often mentioned in the literature about stochastic approximation and recursive stochastic algorithms, this time mainly in its decreasing step version (see [5]), as the *two-armed bandit*. In fact, from a mathematical point of view, it is one of the simplest examples of a stochastic approximation algorithm having a "noiseless trap." We will come back further on this property which was another motivation for investigating this algorithm.

The sequence $(U_n)_{n\geq 1}$ and the events $A_n, B_n, n \geq 1$ are defined on a probability space $(\Omega, \mathcal{A}, \mathbb{P})$. We will make some further assumptions on the events $A_n$ and $B_n$, namely that the sequence

$$(\mathbf{1}_{A_n}, \mathbf{1}_{B_n})_{n\geq 1} \quad \text{is i.i.d.}$$

This assumption corresponds to a "stationary" situation: the traders' daily performances are supposed to be independent and "statistically invariant," that is, identically distributed: so one sets

$$\mathbb{P}(A_1) = p_A \quad \text{and} \quad \mathbb{P}(B_1) = p_B.$$

Of course, the owners of the fund do not know whether $p_A > p_B$ or $p_A < p_B$. Finally, one assumes that the sequences

$$(U_n)_{n\geq 1} \quad \text{and} \quad (\mathbf{1}_{A_n}, \mathbf{1}_{B_n})_{n\geq 1} \quad \text{are independent,}$$

that is, the daily tosses are in no way influenced by the respective past (and future) performances of $A$ and $B$ except for the shares respectively managed that day.

To elucidate the a.s. asymptotic behavior of this allocation procedure, one could call upon classical stochastic approximation methods like the so-called ordinary differential equation (ODE) method. It consists in comparing the asymptotic behavior of the algorithm $(X_n)_{n\geq 1}$ with that of the related $ODE \equiv \dot{x} = \pi h(x)$ where $\pi := p_A - p_B$ and $\pi h(x) := \frac{1}{\gamma_{n+1}} \mathbb{E}(X_{n+1} - X_n | X_n = x) = \pi x(1-x)$ is the *mean* function of the algorithm (see Section 2). One readily checks that this ODE admits two equilibria, 0 and 1, and that, when $p_A > p_B$, its flow $\Phi(x,t) = \frac{x}{(1-x)e^{-\pi t}+x}$ uniformly converges on compact sets of $(0,1]$ toward 1 as $t \to \infty$: the equilibrium 1 is stable with an attraction interval $(0,1]$; on the other hand 0 is repulsive (0 is then called a *trap* for the algorithm). Thus, the celebrated "conditional convergence" theorem due to Kushner and Clark in [8] says that, under technical assumptions fulfilled here, almost every path of the algorithm that visits infinitely many times a compact subset of the attracting interval of a stable equilibrium will converge toward it. Applying that to a path of the two-armed bandit algorithm shows that *if it does not converge to* 0, then it necessarily visits infinitely often the compact interval $[\varepsilon, 1]$ for some $\varepsilon > 0$ and, hence, converges toward 1.

In some way it is not really surprising that this approach fails since stability is a second-order property, whereas the ODE method is based on a first-order approximation. Recent sophisticated first-order approaches like [2] cannot be more efficient for the same reason.

There is a wide literature in stochastic approximation about traps and how not to fall into them (see [3, 6, 9, 14, 16, 17]). They all rely on the fact that, if the noise is exciting enough at a repulsive equilibrium $x^*$, then a.s., the algorithm will not converge to it. By "exciting enough" one means that a conditional variance term at $x^*$ is positive. But the main feature of the two-armed bandit algorithm is that its two equilibria (0 and 1) lie at the boundary of its state space $[0,1]$, so the above conditional variance term is necessarily identically 0 at the repulsive equilibrium $x^* = 0$ (and at $x^* = 1$ as well). So, the behavior of the two-armed bandit algorithm cannot be solved using these approaches.

As far as we know, from a mathematical point of view, the asymptotic behavior of the algorithm has not been elucidated in the literature. The present paper derives from results obtained independently by the third author in [17] and the other two authors.

Heuristics, probably suggested by the behavior of the mean algorithm, seems to consider that the procedure described above works well in practice.

It is interesting for both theoretical and practical motivations to analyze the behavior of the two-armed bandit algorithm, that is:

• Is it possible to choose the gain parameter sequence so that the algorithm a.s. never fails?

• Conversely, does the algorithm "fall in its noiseless trap 0" for some sequences of gain parameters?

This leads to introduce the following terminology when $0 \leq p_B < p_A \leq 1$. (Inverting the rôle played by $A$ and $B$ solves the case $0 \leq p_A < p_B \leq 1$.) The two-armed bandit algorithm is:

• *fallible* when starting from $x \in (0,1)$ if $\mathbb{P}_x(X_n \to 0) > 0$,

- a.s. *infallible* if $\mathbb{P}_x(X_n \to 0) = 0$ for every $x \in (0,1)$.

Although not directly interested by the critical case $p_A = p_B$, we will deeply investigate it since it is a key to solve the general case thanks to a comparison result.

The paper is organized as follows. In Section 1 is stated the main theoretical result of the paper, namely Theorem 1, concerning the convergence and the fallibility of the algorithm. Two corollaries show its consequences on usual parametrized families of steps for which some necessary and sufficient conditions of infallibility are derived.

Section 2 is devoted to some elementary, although important, facts on which relies the proof of Theorem 1, Proposition 2 on one hand and the comparison result stated in Proposition 3 on the other hand. Section 3 is mainly devoted to the proof of items (b) and (c) of Theorem 1 [item (a) is elementary]: Section 3.1 solves item (b) and Section 3.3 solves item (c). Section 3.2 has a particular status: it is a kind of bridge between Sections 3.1 and 3.3: we focus on the special case where the step $\gamma_n$ is constant which is the historical setting considered by those who devised the procedure. It is shown in Theorem 2 that the (positive) probability of failure for the algorithm with constant step $\gamma$ goes to 0 as $\gamma$ goes to zero. Some bounds are displayed, the optimality of which are not known to us. Section 4 makes a connection between regular Pólya urns and the two-armed bandit algorithm: we show that the two-armed bandit algorithm can be seen as a generalized Pólya urn. Thus, we retrieve partially the infallibility results of Theorem 1 using standard methods of proof for the Pólya urns like the "moment method" and the log-method. In the martingale case ($p_A = p_B$) these approaches yield some more information about the distribution of the a.s. limit $X_\infty$ of $X_n$. In Section 5 some first elements about the rate of convergence of the algorithm are provided that emphasize its nonstandard behavior among stochastic approximation procedures. Furthermore, some stopping rules are derived for the algorithm, inspired by some method of proof for infallibility. The last section contains some provisional remarks and additional results.

Note that, except for the notations and the elementary facts contained in Section 2, other sections are self-contained and can be read independently.

NOTATION. (i) The letter $C$ will denote a positive real constant that may change from line to line.

(ii) The letter $\xi$ will denote a *random* positive real constant that may change from line to line.

(iii) Let $(a_n)_{n \geq 0}$ and $(b_n)_{n \geq 0}$ be two sequences of positive real numbers. The symbol $a_n \asymp b_n$ is for $a_n = O(b_n)$ and $b_n = O(a_n)$, whereas the symbol $a_n \sim b_n$ means $\lim_n a_n/b_n = 1$.

**1. The main result.**

THEOREM 1. (a) *Almost sure convergence.*

(i) *If $0 < p_B < p_A \leq 1$ and $x \in (0,1)$, $(X_n)_{n \geq 0}$ is a bounded submartingale, hence $\mathbb{P}_x$-a.s. converging toward a random variable $X_\infty$. The random variable $X_\infty$ takes values in $\{0,1\}$ and*

$$\mathbb{P}_x(X_\infty = 1) = x + \pi \sum_{n \geq 0} \gamma_{n+1} \mathbb{E}_x(h(X_n)) > x + \pi \gamma_1 x(1-x) > x.$$

*(If $p_B = 0$ and $p_A > 0$, then, $X_n$ is nondecreasing and converges toward 1.)*

(ii) *If $0 < p_B = p_A \leq 1$ and $x \in (0,1)$, then $(X_n)_{n \geq 0}$ is a bounded martingale $\mathbb{P}_x$-a.s. converging toward a random variable $X_\infty$.*

*Moreover, if $\sum_{n \geq 0} \gamma_{n+1}^2 = +\infty$, $X_\infty$ is $\{0,1\}$-valued with distribution Bernoulli($x$). (If $p_B = p_A = 0$, then $X_n = x$, $\mathbb{P}_x$-a.s. for every $n \geq 0$.)*

(iii) *If $x \in \{0,1\}$, then $X_n = x$, $\mathbb{P}_x$-a.s. for every $n \geq 0$.*

(b) *Convergence to 0 with positive probability.* *If*

$$\sum_{n \geq 0} \prod_{k=1}^{n} (1 - p_B \gamma_k) < +\infty \tag{3}$$

*then, for every $x \in [0,1)$*

$$\mathbb{P}_x(X_\infty = 0) > 0.$$

*In particular:*

(i) *if $0 < p_B < p_A \leq 1$, then, for every $x \in (0,1)$, the two-armed bandit algorithm starting from $x$ is fallible;*

(ii) *if $0 < p_B = p_A \leq 1$ and $\sum_n \gamma_n^2 < +\infty$, then, for every $x \in (0,1)$,*

$$\mathbb{P}_x(X_\infty = 0), \quad \mathbb{P}_x(X_\infty = 1) \quad \text{and} \quad \mathbb{P}_x(X_\infty \in (0,1)) > 0.$$

(c) *Convergence to a nonzero value.* *Assume $0 \leq p_B \leq p_A \leq 1$ and*

$$\gamma_n = O(\Gamma_n e^{-p_B \Gamma_n}). \tag{4}$$

*Then, for every $x \in (0,1]$,*

$$\mathbb{P}_x(X_\infty = 0) = 0.$$

*In particular:*

(i) *if $0 \leq p_B < p_A \leq 1$ then, for every $x \in (0,1]$,*

$$X_\infty = 1 \quad \mathbb{P}_x\text{-a.s., that is, the algorithm is a.s. infallible,}$$

(ii) *when $0 \leq p_B = p_A \leq 1$ then, for every $x \in (0,1)$,*

$$X_\infty \in (0,1), \quad \mathbb{P}_x\text{-a.s.}$$

PROOF. This theorem follows from Propositions 2, 4 and 5. These propositions can be seen as steps of the proofs of the theorem. □

We will derive in Corollaries 1 and 2 how the above step assumptions (3) and (4) for fallibility or infallibility read on some natural parametrized families of step sequences.

But, first, we will shortly enlighten some connections between the different step assumptions appearing in the statements of the above Theorem 1.

(i) $\gamma_n = O(\Gamma_n e^{-p_B \Gamma_n}) \Rightarrow \sum_n \gamma_n^2 < +\infty$: see the remark after Proposition 5, Section 3.3.

(ii) $\sum_n \prod_{k=1}^n (1 - p_B \gamma_k) = +\infty \not\Rightarrow \sum_n \gamma_n^2 < +\infty$: a counter-example is provided by

$$\gamma_{2^n} = \frac{1}{\sqrt{n+1}}, \qquad n \geq 0 \text{ and } \gamma_k = 0 \text{ if } k \notin \{2^n, n \geq 0\}.$$

COROLLARY 1 (Fallibility). *Let $p_B \in (0,1]$.*

(a) *Constant step. If the step $\gamma_n := \gamma \in (0,1)$, the two-armed bandit algorithm does converge toward $0$ with positive probability. Namely,*

$$\forall x \in (0,1), \qquad \mathbb{P}_x(X_\infty = 0) \geq (1-x)^{1/p_B \gamma} > 0.$$

(b) *Power step* (I). *One considers the family of "power" steps $\gamma_n := (\frac{C}{n+C})^\alpha$, $0 < \alpha \leq 1$, $C > 0$, $n \geq 1$. These step sequences satisfy assumption* (2). *If*

$$(0 < \alpha < 1) \quad or \quad (\alpha = 1 \text{ and } C > 1/p_B),$$

*then, for every $x \in [0,1)$, $\mathbb{P}_x(X_\infty = 0) > 0$ (i.e., the algorithm is fallible from $x$).*

(c) *In particular, if $0 < p_B < p_A \leq 1$, the two-armed bandit algorithm is fallible starting from any $x \in [0,1)$ for the step sequences specified in the above items* (a) *and* (b).

PROOF. (b) The above condition on $C$ and $\alpha$ implies that assumption (3) of Theorem 1 is fulfilled.

(a) The lower bound for $\mathbb{P}_x(X_\infty = 0)$ needs further care. It relies on (9) established in the proof of Proposition 4: setting $\gamma_n = \gamma \in (0,1)$, it reads

$$\mathbb{P}_x(X_\infty = 0) \geq \mathbb{E}_x\left(\prod_{n \geq 1}\left(1 - x \prod_{k=1}^{n-1}(1 - \mathbf{1}_{B_k} \gamma_k)\right)\right).$$

Then the computations can easily be carried on: the Jensen inequality yields

$$\mathbb{P}_x(X_\infty = 0) \geq \exp\left(\sum_{n \geq 1} \mathbb{E}_x \log(1 - xZ_n)\right) \quad \text{with } Z_n := \prod_{k=1}^{n-1}(1 - \gamma \mathbf{1}_{B_k}).$$

Now

$$\mathbb{E}_x(\log(1 - xZ_n)) = -\sum_{m\geq 1} \frac{x^m}{m}\mathbb{E}(Z_n^m)$$

$$= -\sum_{m\geq 1} \frac{x^m}{m}((1-\gamma)^m p_B + 1 - p_B)^{n-1}$$

so that

$$\sum_{n\geq 1}\mathbb{E}_x(\log(1-xZ_n)) = -\frac{1}{p_B}\sum_{m\geq 1}\frac{x^m}{m(1-(1-\gamma)^m)}$$

$$\geq \frac{1}{p_B\gamma}\log(1-x).$$

Finally, for every $x \in [0,1]$,

(5) $$\mathbb{P}_x(X_\infty = 0) \geq (1-x)^{1/p_B\gamma}. \qquad \square$$

COROLLARY 2 (Infallibility). *Let* $0 < p_B \leq p_A \leq 1$.

(a) *Power step* (II). *Let* $\gamma_n := (\frac{C}{n+C})^\alpha$, $0 < \alpha \leq 1$, $C > 0$, $n \geq 1$. *Then,* $\mathbb{P}_x(X_\infty = 0) = 0$ *for every* $x \in (0,1]$ *if and only if*

$$\alpha = 1 \quad \text{and} \quad C \leq \frac{1}{p_B}.$$

(b) *Power step* (III). *Set* $\gamma_n := \frac{\Delta_n}{1+\Delta_1+\cdots+\Delta_n}$, $n \geq 1$, *where* $(\Delta_n)_{n\geq 1}$ *is a sequence of positive real numbers satisfying* $\Delta_n \sim Cn^{1/p_B - 1}\log^\alpha n$ *for some* $\alpha > 0$ *and* $C > 0$. *These step sequences satisfy assumption* (2) *since* $\gamma_n \sim \frac{1}{np_B}$. *Then,* $\mathbb{P}_x(X_\infty = 0) = 0$ *for every* $x \in (0,1]$ *if and only if*

$$\alpha \leq \frac{1}{p_B}.$$

(c) *In particular, if* $0 < p_B < p_A \leq 1$, *the two-armed bandit algorithm is a.s. infallible for the step sequences specified in the above items* (a) *and* (b), *that is,*

$$\forall x \in (0,1], \qquad \mathbb{P}_x(X_\infty = 1) = 1.$$

*Note for practical implementation that the step sequence* $\gamma_n = \frac{1}{n+1}$, $n \geq 1$, *corresponding to* $C = 1$ *always satisfies item* (a) *regardless of the value of* $p_B$ *since* $1 < 1/p_B$.

PROOF. (a) First, assumption (2) is clearly fulfilled. Now, in view of Corollary 1(b), we just need to prove that assumption (4) of Theorem 1 is satisfied if $\alpha = 1$ and $Cp_B \leq 1$. If $\alpha = 1$, $\Gamma_n = C\log n + C' + o(1)$. Consequently, assumption (4) reads $1/n = O(\log(n)n^{-Cp_B})$, that is, $Cp_B \leq 1$.

(b) One has $S_n := 1 + \Delta_1 + \cdots + \Delta_n \sim Cp_B n^{1/p_B} \log^\alpha n$ so that $\gamma_n \sim \frac{1}{p_B} \frac{1}{n}$. Now assume that $\alpha p_B \leq 1$. We need to check assumption (4) of Theorem 1. Notice that $S_n \asymp e^{\Gamma_n}$ [see the preliminary remark after Lemma 1 in Section 3.3 for more details, especially (15)]. Therefore, assumption (4) reduces to $\gamma_n = O(\Gamma_n S_n^{-p_B})$, which follows from $\Gamma_n \sim \frac{1}{p_B} \log n$ and $S_n^{-p_B} \sim (Cp_B)^{-p_B} \frac{1}{\log^{\alpha p_B} n}$. We now prove that if $\alpha p_B > 1$, assumption (3) holds. We have

$$\prod_{k=1}^{n}(1 - p_B \gamma_k) \leq e^{-p_B \Gamma_n} = O(S_n^{-p_B}) = O\left(\frac{1}{n \log^{\alpha p_B} n}\right). \qquad \square$$

FURTHER REMARKS ON THE STEP ASSUMPTIONS. (i) It follows from Corollary 2 that there exist sequences of steps $\gamma_n$ and $\gamma_n'$ satisfying $\gamma_n \sim \gamma_n' \sim \frac{1}{p_B n}$ and such that the corresponding algorithms $X_n$ and $X_n'$ are fallible and infallible, respectively. In fact, the critical case for infallibility is not entirely elucidated by the above results.

(ii) The asymptotics of the constant step setting when the step $\gamma$ goes to 0 is elucidated in Theorem 2, Section 3.2.

**2. Some elementary facts.** The random innovation at time $n$ is clearly $\varepsilon_n := (U_n, \mathbf{1}_{A_n}, \mathbf{1}_{B_n})$ (the $\varepsilon_n$'s are i.i.d.). Set $\mathcal{F}_n := \sigma(\varepsilon_1, \ldots, \varepsilon_n)$, $n \geq 1$ and $\mathcal{F}_0 := \{\varnothing, \Omega\}$. We denote by $\underline{\mathcal{F}}$ the filtration $(\mathcal{F}_n)_{n \geq 0}$. It follows from (1) that $(X_n)_{n \geq 1}$ is obviously a $[0,1]$-valued $\underline{\mathcal{F}}$-Markov chain (homogeneous if $\gamma_n = \gamma$). For notational convenience, we will denote by $\mathbb{P}_x$ the distribution of the whole sequence $(X_n)_{n \geq 0}$ starting at $x \in [0,1]$. One also derives from (1) some straightforward properties of the algorithm.

PROPOSITION 1. *For every $x \in (0,1)$ and every $n \geq 1$, $X_n \in (0,1)$. On the other hand, both states 0 and 1 are absorbing, that is, if $x \in \{0,1\}$, $X_n = x$ for every $n \geq 1$, $\mathbb{P}_x$-a.s.*

Viewed as a stochastic approximation procedure, its canonical form reads

(6) $$X_n = X_{n-1} + \gamma_n \pi h(X_{n-1}) + \gamma_n \Delta M_n,$$

where

$$\pi := p_A - p_B, \qquad h(x) := x(1-x)$$

and

$$\Delta M_n := (1 - X_{n-1})\mathbf{1}_{\{U_n \leq X_{n-1}\} \cap A_n} - X_{n-1}\mathbf{1}_{\{U_n > X_{n-1}\} \cap B_n} - \pi h(X_{n-1})$$

is an $\underline{\mathcal{F}}$-martingale increment.

REMARK. The *mean* algorithm associated with $(X_n)_{n\geq 1}$ is the deterministic recursive procedure defined by

$$x_{n+1} = x_n + \gamma_{n+1}\pi h(x_n), \qquad x_0 \in [0,1].$$

It can be solved very easily: when $\pi = p_A - p_B > 0$ and $x_0 \in (0,1]$, the sequence $x_n$ is $[0,1]$-valued and nondecreasing, hence converging toward $x_\infty$. Since the series $\sum_n \gamma_n h(x_{n-1}) < +\infty$, whereas $\sum_n \gamma_n = +\infty$, it is obvious that $h(x_\infty) = 0$. Hence, $x_\infty = 1$ since $x_\infty \geq x_0 > 0$. So, the mean algorithm never fails in pointing out the best trader since it asymptotically assigns the whole fund to be managed by $A$ when $p_A > p_B$ (and by $B$ when $p_B > p_A$). Unfortunately, it needs to know *a priori* who is the best trader, that is, whether $p_A > p_B$ or $p_B > p_A$.

Similarly, (6) shows that $(X_n)_{n\geq 1}$ is a bounded submartingale and one derives (see Proposition 2) that then $\mathbb{P}_x$-a.s. $X_n$ converges toward a $\{0,1\}$-valued random variable $X_\infty$ if $p_A \neq p_B$. But this time, there is no straightforward argument showing that the procedure always points out the best trader, for example, $X_\infty = 1$ $\mathbb{P}_x$-a.s. when $p_A > p_B$ and $x \in (0,1]$. The next proposition yields some first answers about the behavior of the algorithm.

PROPOSITION 2. (a) *Submartingale case.* If $0 < p_B < p_A \leq 1$ and $x \in (0,1)$, $(X_n)_{n\geq 0}$ is a bounded $\underline{\mathcal{F}}$-submartingale, hence $\mathbb{P}_x$-a.s. converging toward a random variable $X_\infty$, taking values in $\{0,1\}$ and

$$\mathbb{P}_x(X_n \to 1) = x + \pi \sum_{n\geq 0} \gamma_{n+1}\mathbb{E}_x(h(X_n)) > x + \pi\gamma_1 x(1-x) > x.$$

*(If $p_B = 0$ and $p_A > 0$, then, $X_n$ is nondecreasing and converges toward 1.)*

(b) *Martingale case.* If $0 < p_B = p_A \leq 1$ and $x \in (0,1)$, then $(X_n)_{n\geq 0}$ is a bounded $\underline{\mathcal{F}}$-martingale $\mathbb{P}_x$-a.s. converging toward a random variable $X_\infty$. Moreover:

   (i) *if $\sum_{n\geq 0}\gamma_{n+1}^2 = +\infty$, then $X_\infty$ is $\{0,1\}$-valued with Bernoulli distribution $\mathcal{B}(x)$,*

   (ii) *if $\sum_{n\geq 0}\gamma_{n+1}^2 < +\infty$, then $X_\infty$ is $[0,1]$-valued and satisfies $\mathbb{P}_x(X_\infty \in (0,1)) > 0$.*

*(If $p_B = p_A = 0$, then $X_n = x$, $\mathbb{P}_x$-a.s.)*

PROOF. (a) $(X_n)_{n\geq 0}$ is obviously a bounded $\underline{\mathcal{F}}$-submartingale. Furthermore, its a.s. limit, say $X_\infty$, satisfies

$$\mathbb{E}_x(X_\infty) = \lim_n \mathbb{E}_x(X_n)$$

$$= x + \lim_n \mathbb{E}_x\left(\sum_{k=0}^{n-1} \gamma_{k+1}\pi h(X_k) + \gamma_{k+1}\Delta M_{k+1}\right)$$

$$= x + \pi \sum_{n \geq 0} \gamma_{n+1} \mathbb{E}_x h(X_n).$$

Hence, $\sum_{n \geq 0} \gamma_{n+1} \mathbb{E}_x(h(X_n)) < +\infty$ since $\pi > 0$ and consequently, $\sum_{n \geq 0} \gamma_{n+1} h(X_n) < +\infty$ a.s., which in turn implies $\liminf_n h(X_n) = 0$ since $h$ is nonnegative and $\sum_{n \geq 0} \gamma_{n+1} = +\infty$. It follows that $h(X_\infty) = 0$ so that $X_\infty \in \{0, 1\}$ a.s. Finally, $\mathbb{P}_x(X_\infty = 1) = \mathbb{E}_x(X_\infty)$.

(b) $(X_n)_{n \geq 0}$ is obviously a bounded $\underline{\mathcal{F}}$-martingale. When $p_A = p_B$, an elementary computation shows that

$$\mathbb{E}_x(X_n(1 - X_n)) = (1 - p_A \gamma_n^2) \mathbb{E}_x(X_{n-1}(1 - X_{n-1}))$$
(7)
$$= x(1 - x) \prod_{k=1}^{n} (1 - p_A \gamma_k^2)$$

so that

$$\mathbb{E}_x(X_\infty(1 - X_\infty)) = x(1 - x) \prod_{n \geq 1} (1 - p_A \gamma_n^2).$$

The announced result follows since the infinite product converges toward a nonzero limit iff $\sum_n \gamma_n^2 < +\infty$. □

One may specify without loss of generality the definition of events $A_n$ and $B_n$: these two events never interact so only the marginal distributions of $\mathbf{1}_{A_1}$ and $\mathbf{1}_{B_1}$ are involved in the distribution of the whole sequence $(X_n)_{n \geq 0}$. So, one sets

(8) $$A_n := \{V_n \leq p_A\} \quad \text{and} \quad B_n := \{V_n \leq p_B\},$$

where $(U_n)_{n \geq 1}$, $(V_n)_{n \geq 1}$ are two independent i.i.d. $U([0, 1])$-distributed sequences.

Then, this "coupled" algorithm is pathwise monotonous as a function of $p_A$, the parameter $p_B$ being fixed. This is established in the proposition below.

PROPOSITION 3 (Pathwise comparison result). *Let $x \in (0, 1)$. Let $(X_n)$ and $(X'_n)$ denote two "coupled" two-armed bandit algorithms built from the sequences $(U_n)$ and $(V_n)$, starting from $x \leq x'$ and associated to the parameters $(p_B, p_A)$ and $(p_B, p'_A)$, respectively, with $p_A \leq p'_A$. Then for every $n \in \mathbb{N}$,*

$$X_n \leq X'_n.$$

*In particular,*

$$\{X'_n \to 0\} \subset \{X_n \to 0\}.$$

PROOF. The result follows from what happens between time 0 and 1. One inspects the four possible cases following:

(i) On $\{U_1 \leq x'\} \cap \{V_1 \leq p'_A\}$, $X'_1 = x' + \gamma_1(1-x')$ and $X_1 \leq x + \gamma_1(1-x) \leq X'_1$.
(ii) On $\{U_1 \leq x'\} \cap \{V_1 > p'_A\}$, $X'_1 = x'$ and $X_1 \leq x \leq x'$.
(iii) On $\{U_1 > x'\} \cap \{V_1 \leq p_B\}$, $X'_1 = (1-\gamma_1)x'$ and $X_1 = (1-\gamma_1)x \leq X'_1$.
(iv) On $\{U_1 > x'\} \cap \{V_1 > p_B\}$, $X'_1 = x' \geq x = X_1$. □

REMARK. One checks that, when $p_A \geq p_B$, the trajectories of the general form of the algorithm are nondecreasing as a function of their starting value. In particular, *the function* $x \mapsto \mathbb{P}_x(X_\infty = 1)$ *is nondecreasing.*

## 3. When does the two-armed bandit algorithm fail?

3.1. *Quite often...*

PROPOSITION 4. *If* $\sum_{n \geq 0} \prod_{k=1}^n (1 - p_B \gamma_k) < +\infty$ *then,*

$$\mathbb{P}_x(X_n \to 0) > 0 \quad \text{for every } x \in [0,1).$$

PROOF. One considers the event

$$D_\infty := \left\{ U_n > x \prod_{k=1}^{n-1} (1 - \gamma_k \mathbf{1}_{B_k}) \text{ for every } n \geq 1 \right\},$$

where $\prod_{k=1}^0 = 1$. One checks by induction that, on $D_\infty$, $X_n = x \prod_{k=1}^n (1 - \gamma_k \mathbf{1}_{B_k})$. The algorithm is nonincreasing toward its limit $X_\infty$. Hence,

$$\mathbb{E}_x(\mathbf{1}_{D_\infty} X_\infty) \leq \lim_n \mathbb{E}_x(\mathbf{1}_{D_\infty} X_n) \leq \lim_n \mathbb{E}_x\left(x \prod_{k=1}^n (1 - \gamma_k \mathbf{1}_{B_k})\right)$$

$$= x \prod_{n \geq 1} (1 - p_B \gamma_n) = 0$$

since $p_B > 0$, so that $X_\infty = 0$ sur $D_\infty$. On the other hand, $(U_n)_{n \geq 1}$ being i.i.d., uniformly distributed and independent of the sequence $(B_n)_{n \geq 1}$,

$$\mathbb{P}_x(D_\infty) = \mathbb{E}_x(\mathbb{P}_x(D_\infty / \sigma(B_n, n \geq 1))) = \mathbb{E}_x\left(\prod_{n \geq 1}\left(1 - x \prod_{k=1}^{n-1}(1 - \mathbf{1}_{B_k} \gamma_k)\right)\right).$$

Consequently,

(9) $$\mathbb{P}_x(X_\infty = 0) \geq \mathbb{E}_x\left(\prod_{n \geq 1}\left(1 - x \prod_{k=1}^{n-1}(1 - \mathbf{1}_{B_k} \gamma_k)\right)\right).$$

Now, the events $B_k$ are independent, hence

$$\mathbb{E}_x\left(\sum_{n\geq 1}\prod_{k=1}^{n-1}(1-\mathbf{1}_{B_k}\gamma_k)\right)=\sum_{n\geq 1}\prod_{k=1}^{n-1}(1-p_B\gamma_k)<+\infty$$

so that $\sum_{n\geq 1}\prod_{k=1}^{n-1}(1-\mathbf{1}_{B_k}\gamma_k)<+\infty$ $\mathbb{P}_x$-a.s. Consequently the infinite product $\prod_{n\geq 1}(1-x\prod_{k=1}^{n-1}(1-\mathbf{1}_{B_k}\gamma_k))$ converges toward a $\mathbb{P}_x$-a.s. positive random variable. Hence, $\mathbb{P}_x(D_\infty)>0$. □

3.2. *Especially with constant step although...* The fallibility result obtained in Corollary 1(b) for the algorithm with "slowly" decreasing step is to be compared with the asymptotics of its behavior with constant step. We know from Corollary 1(a) that, if $0<p_B\leq p_A\leq 1$ and $\gamma_n=\gamma\in(0,1)$, then, for every $x\in(0,1]$, the algorithm with step $\gamma$ [denoted $(X_n^\gamma)_{n\geq 0}$ in this paragraph] does fail with positive probability: namely, it converges $\mathbb{P}_x$-a.s. toward a $\{0,1\}$-valued random variable $X_\infty^\gamma$ satisfying $\mathbb{P}(X_\infty^\gamma=0)>0$. The following theorem shows, however, that the probability of failure goes to 0 as $\gamma\to 0$.

The fallibility of the algorithm with constant step and this property is known by specialists in Learning Automata theory (see the discussion in Chapter 5 in [11]), although not clearly established mathematically in full generality.

THEOREM 2. *Assume that $0\leq p_B<p_A\leq 1$ and $\gamma_n=\gamma\in(0,1)$. Then, for every $x\in(0,1]$,*

$$\mathbb{P}_x(X_\infty^\gamma=1)\geq 1-\frac{2p_A\gamma}{\pi(1-\gamma)^2}\left(\frac{1}{x}-1\right),$$

*hence,*

$$\lim_{\gamma\to 0_+}\mathbb{P}_x(X_\infty^\gamma=1)=1.$$

PROOF. We focus on the function

$$\psi_\gamma(x):=\mathbb{E}_x\left(\sum_{n\geq 0}h(X_n)\right)=\sum_{n\geq 0}P_\gamma^n h(x),$$

where $P_\gamma(x,dy)$ denotes the Markov transition probability of the two-armed bandit algorithm with constant step $\gamma\in(0,1)$. A straightforward computation shows that for any function $f:[0,1]\to\mathbb{R}$,

$$P_\gamma(f)(x):=\mathbb{E}_x(f(X_1))$$
$$=p_A x f(x+\gamma(1-x))+p_B(1-x)f(x(1-\gamma))$$
$$+(1-p_A x-p_B(1-x))f(x).$$

At this stage, it is convenient to observe that, for every function $g:[0,1] \to \mathbb{R}$,

$$P_\gamma(gh) = hQ_\gamma(g),$$

where the operator $Q_\gamma$ is defined by

$$Q_\gamma(g)(x) = (1-\gamma)\Big(p_A(x+\gamma(1-x))g(x+\gamma(1-x))$$
$$+ p_B(1-x(1-\gamma))g(x(1-\gamma))\Big)$$
$$+ (1 - p_A x - p_B(1-x))g(x).$$

It is clear that $\psi_\gamma$ satisfies $\psi_\gamma = h\chi_\gamma$, where $\chi_\gamma := \sum_{n\geq 0} Q_\gamma^n(\mathbf{1})$. One shows by successive inductions and a little elementary Calculus that $\chi_\gamma$ and $x \mapsto \frac{1}{\pi\gamma x} - \chi_\gamma(x)$ are absolutely decreasing functions [an infinitely differentiable function $f$ is *absolutely decreasing* on $(0,1)$ if its successive derivatives $f^{(n)}$ satisfy $(-1)^n f^{(n)} \geq 0$ for every $n \geq 0$]. On one hand, one derives that $\psi_\gamma$ is indefinitely differentiable on $(0,1]$ and, on the other hand, that

(10) $\quad 0 \leq \chi_\gamma(x) \leq \dfrac{1}{\pi\gamma x}, \quad |\chi'_\gamma(x)| \leq \dfrac{1}{\pi\gamma x^2} \quad \text{and} \quad 0 \leq \chi''_\gamma(x) \leq \dfrac{2}{\pi\gamma x^3}.$

Then, it follows from the definition of $\psi_\gamma$ that

(11) $\quad \psi_\gamma - P_\gamma \psi_\gamma = h \quad \text{and} \quad \mathbb{P}_x(X_\infty^\gamma = 1) = x + \pi\gamma\psi_\gamma(x).$

The first of these two identities reads

$$p_A x(\psi_\gamma(x) - \psi_\gamma(x + \gamma(1-x))) + p_B(1-x)(\psi_\gamma(x) - \psi_\gamma((1-\gamma)x))$$
$$= x(1-x),$$

that is,

$$p_A x \int_\gamma^0 \psi'_\gamma(x + t(1-x))(1-x)\,dt + p_B(1-x)\int_\gamma^0 \psi'_\gamma(x - tx)(-x)\,dt$$
$$= x(1-x).$$

Hence, simplifying by $x(1-x)$ for every $x \in (0,1)$, yields

$$p_A \int_0^\gamma \psi'_\gamma(x + t(1-x))\,dt - p_B \int_0^\gamma \psi'_\gamma(x - tx)\,dt = -1.$$

Now,

$$\psi'_\gamma(x + t(1-x)) = \psi'_\gamma(x) + \int_0^t \psi''_\gamma(x + s(1-x))(1-x)\,ds,$$

$$\psi'_\gamma(x - tx) = \psi'_\gamma(x) + \int_0^t \psi''_\gamma(x(1-s))(-x)\,ds$$

so that

$$\pi\gamma\psi'_\gamma(x) + p_A(1-x)\int_0^\gamma \int_0^t \psi''_\gamma(x+s(1-x))\,ds\,dt$$

$$+ p_B x \int_0^\gamma \int_0^t \psi''_\gamma(x-sx)\,ds\,dt = -1.$$

Combining a rewriting of this identity with an obvious inequality leads to

$$-\psi'_\gamma(x) = \frac{1}{\pi\gamma}\left(1 + \int_0^\gamma \int_0^t \Big(p_A(1-x)\psi''_\gamma(x+s(1-x))\right.$$

$$\left. + p_B x \psi''_\gamma(x(1-s))\Big)\,ds\,dt\right)$$

$$\geq \frac{1}{\pi\gamma} - \frac{p_A(1-x) + p_B x}{\pi\gamma}\frac{\gamma^2}{2}\sup_{u > x(1-\gamma)}|\psi''_\gamma(u)|.$$

Plugging inequalities (10) in the equality $\psi''_\gamma = \chi''_\gamma h + 2\chi'_\gamma h\prime + \chi_\gamma h''$ yields

$$|\psi''_\gamma(x)| \leq \frac{2(1-x)}{\pi\gamma x^2} + \frac{2|1-2x|}{\pi\gamma x^2} + \frac{2}{\pi\gamma x} \leq \frac{4}{\pi\gamma x^2}$$

and consequently,

$$-\psi'_\gamma(x) \geq \frac{1}{\pi\gamma} - \frac{p_A - \pi x}{\pi\gamma}\frac{2\gamma^2}{\pi\gamma(1-\gamma)^2 x^2}.$$

Now $\psi_\gamma(y) = \psi_\gamma(1) + \int_y^1(-\psi'_\gamma(u))\,du = \int_y^1(-\psi'_\gamma(u))\,du$ ($\psi_\gamma(1) = 0$), hence,

$$\psi_\gamma(x) \geq \frac{1-x}{\pi\gamma} + \frac{2p_A}{\pi^2(1-\gamma)^2}\left(1 - \frac{1}{x}\right).$$

Finally, one comes to

$$\mathbb{P}_\gamma(X_n^\gamma \to 1) \geq x + 1 - x - \frac{2p_A\gamma}{\pi(1-\gamma)^2}\left(\frac{1}{x} - 1\right). \qquad \square$$

3.3. *But not always!*

PROPOSITION 5. *Assume $0 \leq p_B \leq p_A \leq 1$ and assumption* (4). *Then, for every $x \in (0,1]$,*

(12) $$\mathbb{P}_x(X_\infty = 0) = 0.$$

REMARK. As soon as $p_B > 0$, assumption (4) implies $\sum_{n \geq 1} \gamma_n^2 < +\infty$: for $n \geq 1$,

$$\sum_{k=1}^n \gamma_k^2 \leq C \sum_{k=1}^n (\Gamma_k - \Gamma_{k-1})\Gamma_k e^{-p_B \Gamma_k}$$

$$\leq C \int_0^{\Gamma_n} u e^{-p_B u}\,du \leq C \int_0^{+\infty} u e^{-p_B u}\,du < +\infty,$$

since $u \mapsto ue^{-p_B u}$ is nonincreasing for $u$ large enough.

LEMMA 1 (Case $p_A = p_B = 1$). *Assume $p_A = p_B = 1$ and set $\Delta_0 := 1$ and $\Delta_n := \frac{\gamma_n}{\prod_{k=1}^n (1-\gamma_k)}$, for $n \geq 1$. If the sequence $(\Delta_n)_{n \in \mathbb{N}}$ satisfies*

(13) $$\Delta_n = O(\Gamma_n),$$

*then, for every $x \in (0,1]$, $\mathbb{P}_x(X_\infty = 0) = 0$.*

PRELIMINARY REMARK. With the notation of the lemma, we have

(14) $$\gamma_n = \frac{\Delta_n}{S_n} \quad \text{with } S_n = \sum_{k=0}^n \Delta_k.$$

The partial sums $S_n$ and $\Gamma_n$ satisfy, for every $n \geq 1$,

(15) $$\log S_n - \sum_{k=1}^n \frac{\gamma_k^2}{1-\gamma_k} \leq \Gamma_n \leq \log S_n.$$

This follows from the easy comparisons (with $S_0 = \Delta_0 = 1$),

$$\Gamma_n = \sum_{k=1}^n \frac{\Delta_k}{S_k} \begin{cases} \leq \int_1^{S_n} \frac{du}{u} = \log S_n, \\ \geq \sum_{k=1}^n \frac{S_{k-1}}{S_k} \int_{S_{k-1}}^{S_k} \frac{du}{S_{k-1}} \geq \sum_{k=1}^n (1-\gamma_k) \int_{S_{k-1}}^{S_k} \frac{du}{u} \\ \geq \log S_n - \sum_{k=1}^n \frac{\gamma_k^2}{1-\gamma_k}. \end{cases}$$

Hence, $S_n \to \infty$ as $n \to \infty$ since $\Gamma_n \to \infty$ and $\gamma_n \leq C\Gamma_n/S_n \leq C\Gamma_n e^{-\Gamma_n}$. Consequently, the former remark applies with $p_B = 1$ and shows that

$$\sum_{n \geq 1} \gamma_n^2 < +\infty.$$

Finally, assumption (13) implies

(16) $$\Gamma_n \sim \log S_n \quad \text{and} \quad S_n \asymp e^{\Gamma_n} \qquad \text{as } n \to \infty.$$

PROOF OF LEMMA 1. The algorithm can now be rewritten as follows (we assume $p_A = p_B = 1$):

$$S_{n+1} X_{n+1} = S_{n+1} X_n + \Delta_{n+1}(\mathbf{1}_{\{U_{n+1} \leq X_n\}} - X_n).$$

Hence,

$$S_{n+1} X_{n+1} = S_n X_n + \Delta_{n+1} \mathbf{1}_{\{U_{n+1} \leq X_n\}}.$$

Let $Y_n := S_n X_n$. We will first prove that $\lim_n Y_n = +\infty$, a.s. Since the sequence $(\Delta_n/\Gamma_n)_{n\geq 1}$ is bounded and $\mathbb{E}(\mathbf{1}_{\{U_{n+1}\leq X_n\}}|\mathcal{F}_n) = X_n$, we have (see, e.g., Theorem 2.7.33 in [4])

$$\left\{\sum_{n=0}^{\infty} \frac{\Delta_{n+1}}{\Gamma_{n+1}} \mathbf{1}_{\{U_{n+1}\leq X_n\}} = \infty\right\} = \left\{\sum_{n=0}^{\infty} \frac{\Delta_{n+1}}{\Gamma_{n+1}} X_n = \infty\right\} \qquad \text{a.s.}$$

But

$$\Delta_{n+1} X_n \geq \Delta_{n+1} x \prod_{k=1}^{n}(1-\gamma_k) \geq \Delta_{n+1} x \prod_{k=1}^{n+1}(1-\gamma_k) = \gamma_{n+1} x$$

so that

$$\sum_{n=0}^{\infty} \frac{\Delta_{n+1}}{\Gamma_{n+1}} X_n \geq x \sum_{n=0}^{\infty} \frac{\gamma_{n+1}}{\Gamma_{n+1}} = +\infty,$$

since $\Gamma_n \to \infty$ and $\gamma_n \to 0$ as $n \to \infty$. Consequently, the nondecreasing sequence $(Y_n)_{n\geq 0}$ satisfies

$$\limsup_n Y_n = \limsup_n \sum_{k=0}^{n-1} \Delta_{k+1} \mathbf{1}_{\{U_{k+1}\leq X_k\}} \geq \gamma_1 \sum_{n\geq 0} \frac{\Delta_{n+1}}{\Gamma_{n+1}} \mathbf{1}_{\{U_{n+1}\leq X_k\}} = +\infty.$$

Next, we prove that $\limsup_n \frac{Y_n}{\log S_n} = +\infty$ a.s. One may write $Y_n = x + \sum_{k=0}^{n-1} \Delta_{k+1} \mathbf{1}_{\{U_{k+1}\leq Y_k/S_k\}}$ so that for any $\lambda > 0$,

$$\limsup_n \frac{Y_n}{\log S_n} \geq \limsup_n \frac{Z_n^\lambda}{\log S_n} \qquad \text{where } Z_n^\lambda = \sum_{k=0}^{n-1} \Delta_{k+1} \mathbf{1}_{\{U_{k+1}\leq \lambda/S_k\}}.$$

We have

$$\mathbb{E} Z_n^\lambda = \sum_{k=0}^{n-1} \Delta_{k+1} \min\left(1, \frac{\lambda}{S_k}\right)$$

and

$$\mathrm{Var}(Z_n^\lambda) = \sum_{k=0}^{n-1} \Delta_{k+1}^2 \min\left(1, \frac{\lambda}{S_k}\right)\left(1 - \min\left(1, \frac{\lambda}{S_k}\right)\right)$$

$$\leq C \log S_n \sum_{k=0}^{n-1} \Delta_{k+1} \min\left(1, \frac{\lambda}{S_k}\right)$$

$$= C \log S_n \mathbb{E} Z_n^\lambda.$$

Consequently,

$$\mathbb{P}(|Z_n^\lambda - \mathbb{E} Z_n^\lambda| \geq \rho \mathbb{E} Z_n^\lambda) \leq C \frac{\log S_n \mathbb{E} Z_n^\lambda}{\rho^2 (\mathbb{E} Z_n^\lambda)^2} \leq C \frac{\log S_n}{\rho^2 \sum_{k=0}^{n-1} \Delta_{k+1} \min(1, \lambda/S_k)}.$$

One checks that $\lim_n \frac{\sum_{k=0}^{n-1} \Delta_{k+1} \min(1,\lambda/S_k)}{\log S_n} = \lambda$, since $\Delta_n \min(1,\lambda/S_{n-1}) \sim \lambda \frac{\gamma_n}{1-\gamma_n} \sim \lambda \gamma_n$.

Let $A_n^\lambda = \{|Z_n^\lambda - \mathbb{E} Z_n^\lambda| < \rho \mathbb{E} Z_n^\lambda\}$. For $\lambda$ large enough, $\mathbb{P}(A_n^\lambda) \geq 1/2$, so that $\mathbb{P}(\limsup_n A_n^\lambda) \geq 1/2$. Now, on the event $\limsup_n A_n^\lambda$,

$$Z_n^\lambda \geq (1-\rho)\mathbb{E} Z_n^\lambda \geq \lambda(1-\rho)\log S_n \qquad \text{for infinitely many } n.$$

Hence, $\mathbb{P}(\limsup_n \frac{Z_n^\lambda}{\log S_n} \geq \lambda(1-\rho)) \geq 1/2$. But the random variable $\limsup_n \frac{Z_n^\lambda}{\log S_n}$ lies in the asymptotic $\sigma$-field of the i.i.d. random variables $U_n$'s, hence,

$$\limsup_n \frac{Z_n^\lambda}{\log S_n} \geq \lambda(1-\rho), \qquad \mathbb{P}_x\text{-a.s.}$$

This holds for every $\rho > 0$ and $\lambda > 0$ so that $\limsup_n \frac{Y_n}{\log S_n} = +\infty$, $\mathbb{P}_x$-a.s.

On the other hand, for any positive integer $p$,

$$\mathbb{E}((X_\infty - X_p)^2|\mathcal{F}_p) = \mathbb{E}\left(\sum_{k=p}^\infty \gamma_{k+1}^2 X_k(1-X_k)\Big|\mathcal{F}_p\right)$$

$$\leq \mathbb{E}\left(\sum_{k=p}^\infty \gamma_{k+1}^2 X_k\Big|\mathcal{F}_p\right) = X_p \sum_{k=p}^\infty \gamma_{k+1}^2.$$

Now observe that

$$\mathbb{P}(X_\infty = 0|\mathcal{F}_p) = \frac{\mathbb{E}(\mathbf{1}_{\{X_\infty=0\}} X_p^2|\mathcal{F}_p)}{X_p^2} \leq \frac{\mathbb{E}((X_\infty - X_p)^2|\mathcal{F}_p)}{X_p^2}$$

$$\leq \frac{\sum_{k=p}^\infty \gamma_{k+1}^2}{X_p} = S_p \frac{\sum_{k=p}^\infty \gamma_{k+1}^2}{Y_p}$$

$$\leq C\frac{S_p}{Y_p} \sum_{k \geq p+1} \frac{\Gamma_k}{S_k^2} \Delta_k \leq C\frac{S_p}{Y_p} \int_{S_p}^{+\infty} \frac{\log u}{u^2}\,du$$

$$\leq C\frac{S_p}{Y_p} \times \frac{\log S_p}{S_p} = C\frac{\log S_p}{Y_p}.$$

One concludes by noting that the bounded martingale $\mathbb{P}_x(X_\infty = 0|\mathcal{F}_p)$, $\mathbb{P}_x$-a.s., converges to $\mathbf{1}_{\{X_\infty=0\}}$ so that

$$\mathbb{P}(X_\infty = 0) = \lim_p \mathbb{E}(\mathbb{P}(X_\infty = 0|\mathcal{F}_p)) \leq C \liminf_p \frac{\log S_p}{Y_p} = 0. \qquad \square$$

LEMMA 2. *Assume $0 < p_A = p_B < 1$ and assumption* (4). *Then, for every $x \in (0,1]$,*

(17) $$\mathbb{P}_x(X_\infty = 0) = 0.$$

PROOF OF LEMMA 2. In that case, the algorithm can be written as follows [see the specification of the algorithm in (8)]:

$$X_{n+1} = X_n + \gamma_{n+1}\mathbf{1}_{B_{n+1}}(\mathbf{1}_{\{U_{n+1}\leq X_n\}} - X_n) \qquad \text{with } B_{n+1} = \{V_{n+1} \leq p_B\}.$$

By conditioning on the $\sigma$-field generated by the events $B_n$, $n \geq 1$, we easily deduce from Lemma 1 that if the sequences $(\Delta_n^B)_{n\in\mathbb{N}}$, defined by $\Delta_n^B = \gamma_n\mathbf{1}_{B_n}/\prod_{k=1}^n(1-\gamma_k\mathbf{1}_{B_k})$ and $\gamma_n^B = \gamma_n\mathbf{1}_{B_n}$ satisfy assumption (13) in Lemma 1, the announced result is proved.

Now, for $n \geq 1$, define

$$M_n = \sum_{k=1}^n \log(1 - \gamma_k\mathbf{1}_{B_k}) - p_B\log(1-\gamma_k) = \sum_{k=1}^n \log(1-\gamma_k)(\mathbf{1}_{B_k} - p_B).$$

The sequence $(M_n)$ is a martingale and $\sup_n \mathbb{E}M_n^2 < +\infty$ because $\sum_n \gamma_n^2 < \infty$ (as follows from the remark below Proposition 5). Therefore, the ratio

$$\frac{\prod_{k=1}^n(1-\gamma_k)^{p_B}}{\prod_{k=1}^n(1-\gamma_k\mathbf{1}_{B_k})} \qquad \text{is a.s. bounded.}$$

Consequently, there exists a $\sigma(B_n, n \geq 1)$-measurable random positive constant $\xi$ such that, $\mathbb{P}_x$-a.s.,

$$\Delta_n^B \leq \xi\frac{\gamma_n^B}{\prod_{k=1}^n(1-\gamma_k)^{p_B}} = \xi\gamma_n^B S_n^{p_B} \leq \xi\gamma_n S_n^{p_B}.$$

Inequality (15) and $\sum_{n\geq 1}\gamma_n^2 < +\infty$ imply that $S_n \leq Ce^{\Gamma_n}$. In turn, assumption (4) yields

$$\Delta_n^B \leq \xi\Gamma_n e^{-p_B\Gamma_n}(e^{\Gamma_n})^{p_B} \leq \xi\Gamma_n.$$

Now, a straightforward martingale argument shows that, a.s.,

$$\sum_{k=1}^n \gamma_k^B \sim p_B \sum_{k=1}^n \gamma_k$$

so that $\Delta_n^B = O(\gamma_1^B + \cdots + \gamma_n^B)$ and the proposition follows from Lemma 1. □

PROOF OF PROPOSITION 5. Using Proposition 3 (pathwise comparison result), one may assume without loss of generality that $p_A = p_B > 0$. Lemma 2 completes the proof. □

**4. The two-armed bandit algorithm as a generalized Pólya urn.** In this section we propose two proofs of Proposition 5—in some special cases—in which the martingale case is based on methods directly inspired by the Pólya urn.

4.1. *Short background on Pólya urn.* Assume that an urn contains at time 0, $r$ red balls and $b$ black balls. At every time $n$ one draws at random a ball from the urn and then puts it back in the urn with another ball of the same color. Then, at every time $n$, the urn contains (once the new ball has been put in the urn) exactly $r+b+n$ balls. Let $\beta_n$ denote the number of black balls inside the urn at time $n$, let $X_n := \frac{\beta_n}{r+b+n}$ denote the proportion of black balls at time $n$ and $Y_n := \frac{\beta_n}{r+b}$. One models the drawings using a sequence $(U_n)_{n\geq 1}$ of i.i.d. random variables uniformly distributed over $[0,1]$ as follows: if $U_{n+1} \leq X_n$, the ball drawn at time $n+1$ is black, otherwise it is red. Then, these sequences satisfy, respectively,

$$\beta_0 := b \quad \text{and} \quad \beta_{n+1} = \beta_n + \mathbf{1}_{\{U_{n+1} \leq X_n\}},$$

$$Y_0 := \frac{b}{r+b} \quad \text{and} \quad Y_{n+1} = Y_n + \frac{1}{r+b}\mathbf{1}_{\{U_{n+1}\leq X_n\}},$$

$$X_0 := \frac{b}{r+b} \quad \text{and} \quad X_{n+1} = X_n + \frac{1}{r+b+n+1}(\mathbf{1}_{\{U_{n+1}\leq X_n\}} - X_n).$$

Consequently, the regular Pólya urn appears as a special case of the two-armed-bandit algorithm (in the martingale setting $p_A = p_B = 1$) corresponding to a rational starting value $X_0 = \frac{b}{r+b}$ and a step $\gamma_n := \frac{1}{r+b+n}$, that is, $\Delta_n = \frac{1}{r+b}$.

This suggests to try extending some classical methods of proof devised for the Pólya urn to solve the martingale case of the two-armed bandit algorithm, with the hope, in some cases, to get more accurate results, for example, concerning the distribution of the limit $X_\infty$.

4.2. *The moment approach.* Following a classical method devised to solve the Pólya urn (see, e.g., [1]), it is possible to obtain some moment estimates for the limiting distribution of the $X_n$'s in the martingale case $p_A = p_B = 1$. When $\Delta_n = \Delta > 0$, this limiting distribution is even explicit.

PROPOSITION 6. *Assume that $p_A = p_B = 1$ and that the sequence $(\Delta_n)_{n\geq 1}$ is nonincreasing ($\Delta_1$ may be greater than $\Delta_0 = 1$).*

(a) *For every $x \in [0,1]$ and for every integer $m \geq 1$,*

$$\mathbb{E}_x(X_\infty^{m+1}) \leq \prod_{k=0}^{m}\left(1 - \frac{1-x}{S_k}\right) \quad \text{and} \quad \mathbb{E}_x((1-X_\infty)^{m+1}) \leq \prod_{k=0}^{m}\left(1 - \frac{x}{S_k}\right).$$

(18)

*In particular, for every $x \in (0,1)$,*

$$\mathbb{P}_x(X_\infty = 1) \leq \inf_m \mathbb{E}_x(X_\infty^m) = 0 \quad \text{and} \quad \mathbb{P}_x(X_\infty = 0) \leq \inf_m \mathbb{E}_x((1-X_\infty)^m) = 0$$

*since $\sum_{k\geq 1} \frac{1}{S_k} \geq \sum_{k\geq 1} \frac{1}{1+k\Delta_1} = +\infty$.*

(b) *If, moreover* $\Delta_n = \Delta > 0$, $n \geq 1$ *(and $\Delta_0 = 1$), then* $\gamma_n = \frac{\Delta}{n\Delta + 1}$ *and*

$$X_\infty \overset{\mathcal{L}}{\sim} \beta\left(\frac{x}{\Delta}; \frac{1-x}{\Delta}\right).$$

REMARK. Item (b) is a classical result about Pólya urn.

PROOF. (a) One uses the notation $Y_n = S_n X_n$ of Lemma 1 and sets

$$Z_n^{(m)} := \frac{Y_n}{S_n} \frac{Y_n + \Delta_{n+1}}{S_{n+1}} \times \cdots \times \frac{Y_n + \Delta_{n+1} + \cdots + \Delta_{n+m}}{S_{n+m}},$$

$$\mathbb{E}_x(Z_{n+1}^{(m)}/\mathcal{F}_n)$$
$$= \mathbb{E}_x(Z_{n+1}^{(m)}\mathbf{1}_{\{U_{n+1} > X_n\}}/\mathcal{F}_n) + \mathbb{E}_x(Z_{n+1}^{(m)}\mathbf{1}_{\{U_{n+1} \leq X_n\}}/\mathcal{F}_n)$$
$$= \frac{Y_n}{S_{n+1}} \frac{Y_n + \Delta_{n+2}}{S_{n+2}} \times \cdots \times \frac{Y_n + \Delta_{n+2} + \cdots + \Delta_{n+m+1}}{S_{n+m+1}} (1 - X_n)$$
$$+ \frac{Y_n + \Delta_{n+1}}{S_{n+1}} \frac{Y_n + \Delta_{n+1} + \Delta_{n+2}}{S_{n+2}} \times \cdots$$
$$\times \frac{Y_n + \Delta_{n+1} + \cdots + \Delta_{n+m+1}}{S_{n+m+1}} X_n.$$

Recall that $X_n = \frac{Y_n}{S_n}$ and $1 - X_n = \frac{S_n - Y_n}{S_n}$. Hence,

$$\mathbb{E}_x(Z_{n+1}^{(m)}/\mathcal{F}_n)$$
$$= \frac{S_n - Y_n}{S_{n+m+1}} \times \underbrace{\frac{Y_n}{S_n} \frac{Y_n + \Delta_{n+2}}{S_{n+1}} \times \cdots \times \frac{Y_n + \Delta_{n+2} + \cdots + \Delta_{n+m+1}}{S_{n+m}}}_{\leq Z_n^{(m)} \text{ since } \Delta_{n+i+1} \leq \Delta_{n+i}}$$
$$+ \frac{Y_n + \Delta_{n+1} + \cdots + \Delta_{n+m+1}}{S_{n+m+1}}$$
$$\times \underbrace{\frac{Y_n}{S_n} \frac{Y_n + \Delta_{n+1}}{S_{n+1}} \times \cdots \times \frac{Y_n + \Delta_{n+1} + \cdots + \Delta_{n+m}}{S_{n+m}}}_{Z_n^{(m)}}$$
$$\leq Z_n^{(m)}.$$

The sequence $(Z_n^{(m)})_{n \geq 0}$ is then a super-martingale. On the other hand $Z_n^{(m)}$ obviously converges toward $X_\infty^{m+1}$ since, for every $k \leq m$, $\frac{\Delta_{n+k}}{S_{n+m}} \leq \frac{\Delta_1}{S_{n+m}} \to 0$ as $n \to +\infty$ because the sequence $(\Delta_n)_{n \geq 1}$ is nonincreasing and $S_n \uparrow +\infty$ [see (15)] in the preliminary remark following Lemma 1. Consequently, for instance, via Fatou's lemma, for every $m \in \mathbb{N}$,

$$\mathbb{E}_x(X_\infty^{m+1}) \leq Z_0^{(m)} = x \prod_{k=1}^m \frac{x + \Delta_1 + \cdots + \Delta_k}{1 + \Delta_1 + \cdots + \Delta_k}$$

$$= x \prod_{k=1}^m \frac{x - 1 + S_k}{S_k}$$

$$= \prod_{k=0}^m \left(1 - \frac{1-x}{S_k}\right).$$

One proceeds symmetrically with $\widetilde{X}_n := 1 - X_n$ and $\widetilde{Y}_n = S_n - Y_n$ to establish the moment inequalities concerning $1 - X_\infty$.

Hence, the Lebesgue dominated convergence theorem implies that when $m \to \infty$,

$$\mathbb{P}_x(X_\infty = 0) = \lim_m \mathbb{E}_x((1 - X_\infty)^{m+1}) \leq \prod_{k \geq 0} \left(1 - \frac{x}{S_k}\right) = 0$$

since $\sum_{n \geq 0} \frac{1}{S_n} \geq \sum_{n \geq 0} \frac{1}{1 + n\Delta_1} = +\infty$.

(b) When $\Delta_n = \Delta$, $n \geq 1$, the same proof shows that $(Z_n^{(m)})_{n \geq 0}$ is a martingale, hence, for every $m \geq 0$,

$$\mathbb{E}_x(X_\infty^{m+1}) = \prod_{k=0}^m \left(1 - \frac{1-x}{1 + k\Delta}\right) = \prod_{k=0}^m \frac{x/\Delta + k}{(1-x)/\Delta + x/\Delta + k}.$$

Hence, $X_\infty$ has the moments of a $\beta(x/\Delta; (1-x)/\Delta)$ distribution. Both distributions have compact support, hence, they are equal. □

The above result can be extended to the general martingale case.

COROLLARY 3. *If $p_A = p_B \in (0, 1]$ and the sequence $(\Delta_n)_{n \geq 1}$ is nonincreasing ($\Delta_1$ may be greater than $\Delta_0 = 1$). Then, for every $x \in (0, 1)$ $\mathbb{P}_x(X_\infty = 1) = \mathbb{P}_x(X_\infty = 0) = 0$ and, for every $m \in \mathbb{N}$,*

$$(19) \qquad \mathbb{E}_x(X_\infty^{m+1}) \leq \prod_{k=0}^m \left(1 - (1-x) \prod_{\ell=1}^k (1 - \gamma_\ell)\right)$$

$$(20) \qquad \mathbb{E}_x((1 - X_\infty)^{m+1}) \leq \prod_{k=0}^m \left(1 - x \prod_{\ell=1}^k (1 - \gamma_\ell)\right).$$

PROOF. Dealing with the case $p := p_A = p_B < 1$ still needs conditioning with respect to the $\sigma$-algebra generated by the events $B_n$. However, we will proceed slightly differently from in the proof of Proposition 5. We introduce the successive stopping times

$$\forall \omega \in \Omega, \qquad \tau_0^B(\omega) := 0, \qquad \tau_n^B(\omega) := \min\{k > \tau_{n-1}^B / \omega \in B_k\}, \qquad n \geq 1.$$

The $\tau_n$'s are $\mathbb{P}_x$-a.s. finite iff $p_B > 0$. Then, set $\widetilde{\gamma}_n := \gamma_{\tau_n^B}$ and $\widetilde{\Delta}_n$ and $\widetilde{S}_n := 1 + \widetilde{\Delta}_1 + \cdots + \widetilde{\Delta}_n$, $n \geq 1$, as in Lemma 1 so that $\widetilde{\gamma}_n = \widetilde{\Delta}_n / \widetilde{S}_n$. One checks that

$$\frac{\widetilde{\Delta}_{n+1}}{\widetilde{\Delta}_n} = \frac{\Delta_{\tau_{n+1}^B}}{\Delta_{\tau_n^B}} \prod_{k=\tau_n^B+1}^{\tau_{n+1}^B-1} (1 - \gamma_k) \leq \frac{\Delta_{\tau_{n+1}^B}}{\Delta_{\tau_n^B}},$$

so that $\widetilde{\Delta}_n$ is nonincreasing as long as $\Delta_n$ is. It follows from Proposition 6 and obvious equalities that, for every $m \geq 0$,

$$\mathbb{E}_x(X_\infty^{m+1}) \leq \mathbb{E}_x\left(\prod_{k=0}^m \left(1 - \frac{1-x}{\widetilde{S}_k}\right)\right)$$

and

$$\mathbb{E}_x((1 - X_\infty)^{m+1}) \leq \mathbb{E}_x\left(\prod_{k=0}^m \left(1 - \frac{x}{\widetilde{S}_k}\right)\right).$$

Now $\widetilde{S}_n \leq 1 + n\widetilde{\Delta}_1$, $n \geq 1$. Then, proceeding as in the proof of Proposition 6(a) yields

$$\mathbb{P}_x(X_\infty = 0) = \lim_m \mathbb{E}_x((1 - X_\infty)^{m+1}) \leq \prod_{k \geq 0} \left(1 - \frac{x}{1 + k\Delta_1}\right) = 0.$$

One gets similarly that $\mathbb{P}_x(X_\infty = 1) = 0$. The moment bounds follow from the easy fact that $\widetilde{\Delta}_n \leq \Delta_n$ so that $\widetilde{S}_n \leq S_n = (\prod_{k=1}^n (1 - \gamma_k))^{-1}$. □

REMARK. The above bounds (19) and (20) do not involve $p_B$. In fact, they can be improved by replacing $\gamma_n$ by $\gamma_{\tau_n^B}$ in their right-hand side and taking the expectation with respect to $\mathbb{E}_x$. Then, one may use that $\tau_n^B - \tau_{n-1}^B$ is i.i.d. with Geometric distribution $G(p_B)$.

4.3. *The log-martingale approach.* The log-martingale method is another classical approach to Pólya urn. It also yields a new proof of Lemma 1 when the sequence $(\Delta_n)_{n \geq 1}$ is bounded. We use the same notations as in the original lemma.

PROPOSITION 7. *Assume $p_A = p_B = 1$ and $(\Delta_n)_{n \geq 1}$ is bounded.*

(a) *Then, there exists a martingale $(N_n)_{n\geq 1}$ with bounded increments such that, for every $x \in (0,1)$,*

$$\sup_n \left|\log\left(\frac{X_n}{1-X_n}\right) - N_n\right| < +\infty, \qquad \mathbb{P}_x\text{-a.s.}$$

(b) *Consequently, for every $x \in (0,1)$, $\mathbb{P}_x(X_\infty = 0) = \mathbb{P}_x(X_\infty = 1) = 0$.*

PROOF. (a) Set $Z_n := \log(\frac{X_n}{1-X_n}) = \log(\frac{Y_n}{S_n-Y_n}) = \log Y_n - \log(S_n - Y_n)$.

$$\left|\log Y_n - \left(\log x + \sum_{k=1}^n \frac{\Delta Y_k}{Y_{k-1}}\right)\right|$$

$$\leq \sum_{k=1}^n \left|\log\left(\frac{Y_k}{Y_{k-1}}\right) - \frac{\Delta Y_k}{Y_{k-1}}\right| \leq \frac{1}{2}\sum_{k=1}^n \left(\frac{\Delta Y_k}{Y_{k-1}}\right)^2$$

$$\leq \frac{\sup_n \Delta_n}{2} \sum_{k=1}^n \frac{\Delta Y_k}{Y_{k-1}^2} \leq \frac{c^2 \sup_n \Delta_n}{2} \sum_{k\geq 1} \frac{\Delta Y_k}{Y_k^2} < +\infty,$$

where $c := 1 + \sup_n \Delta_n / x$ satisfies $Y_k/Y_{k-1} \leq 1 + \Delta_k/Y_{k-1} \leq c$. Similarly,

$$\left|\log(S_n - Y_n) - \left(\log(1-x) + \sum_{k=1}^n \frac{\Delta_k - \Delta Y_k}{S_{k-1} - Y_{k-1}}\right)\right|$$

$$\leq \frac{\sup_n \Delta_n}{2} \sum_{k\geq 1} \frac{(\Delta_k - \Delta Y_k)}{(S_{k-1} - Y_{k-1})^2} < +\infty.$$

Combining these inequalities yields

$$\sup_n \left|Z_n - \sum_{k=1}^n \left(\frac{\Delta Y_k}{Y_{k-1}} - \frac{\Delta_k - \Delta Y_k}{S_{k-1} - Y_{k-1}}\right)\right| < +\infty, \qquad \mathbb{P}_x\text{-a.s.}$$

Now $N_n := \sum_{k=1}^n (\frac{\Delta Y_k}{Y_{k-1}} - \frac{\Delta_k - \Delta Y_k}{S_{k-1} - Y_{k-1}})$ is a martingale since

$$\mathbb{E}_x(\Delta N_n / \mathcal{F}_{n-1})$$

$$= \mathbb{E}_x(\Delta Y_n / \mathcal{F}_{n-1})\frac{1}{Y_{n-1}} - (\Delta_n - \mathbb{E}_x(\Delta Y_n / \mathcal{F}_{n-1}))\frac{1}{S_{n-1} - Y_{n-1}}$$

$$= \frac{X_{n-1}\Delta_n}{Y_{n-1}} - \frac{(1-X_{n-1})\Delta_n}{S_{n-1} - Y_{n-1}} = \frac{\Delta_n}{S_{n-1}} - \frac{\Delta_n}{S_{n-1}} = 0.$$

Furthermore, its increments are bounded. As a matter of fact $\Delta Y_k$ and $\Delta_k - \Delta Y_k$ are never simultaneously zero and are upper-bounded by $\Delta_k$; on the other hand, $Y_{k-1}$ and $S_{k-1} - Y_{k-1}$ are lower bounded by $x$ and $1-x$, respectively, hence, for every $k \geq 1$,

$$|\Delta N_k| \leq \frac{\Delta_k}{Y_{k-1} \wedge (S_{k-1} - Y_{k-1})} \leq \frac{\|\Delta\|_\infty}{x \wedge (1-x)}.$$

(b) Let $(\langle N\rangle_n)_{n\geq 1}$ denote the conditional variance increment process of the martingale $(N_n)_{n\geq 1}$ and let $\langle N\rangle_\infty$ denote its limit as $n$ goes to infinity. The law of iterated logarithm for martingales with bounded increments says that, on the event $\{\langle N\rangle_\infty = +\infty\}$, the martingale $(N_n)$ satisfies $\liminf_n N_n = -\infty$ and $\limsup_n N_n = +\infty$ a.s. Meanwhile $Z_n$ converges toward $\log(\frac{X_\infty}{1-X_\infty}) \in \overline{\mathbb{R}}$ a.s. The difference of these two quantities remaining bounded, it follows that $\langle N\rangle_\infty < +\infty$, $\mathbb{P}_x$-a.s. Hence, the martingale $N_n$ converges toward a finite limit and, consequently, $X_\infty \in (0,1)$, $\mathbb{P}_x$-a.s. $\square$

REMARK. The above assumption in Proposition 7 does not embody the Power step (III) setting of Corollary 2(b), ($p_A = p_B = 1$ and) $\Delta_n \sim C \log n$, that is, the closest case to the critical case that we can get.

The extension to general $p_A$ and $p_B$, $0 < p_B \leq p_A$, in that framework consists in proving that the sequence $(\Delta_n^B)_{n\geq 1}$ is a.s. bounded. One shows using martingale methods of Lemma 2 that this leads to the condition $\gamma_n = O(e^{-p_B \Gamma_n})$ which is, as expected, more stringent than assumption (4).

## 5. Rate of convergence, stopping rules.

5.1. *Rate of convergence.* The aim of this section is not to elucidate completely the rate of convergence of the two-armed bandit algorithm but to draw some first conclusions from some by-products of the convergence proof. They emphasize that the two-armed bandit algorithm does not behave like a standard stochastic approximation algorithm in terms of rate of convergence. In particular, in some natural situations it may converge infinitely faster than its associated deterministic algorithm in average. This enlightens that the usual CLT for stochastic algorithms proposed in the literature (see, e.g., [5]) does not apply.

First, let us have a look at the algorithm *in average*,

$$x_{n+1} = x_n + \pi \gamma_{n+1} x_n (1 - x_n), \qquad x_0 = x \in (0,1) \qquad \text{with } \pi = p_A - p_B > 0.$$

One has by a straightforward induction that the sequence $(x_n)_{n\geq 0}$ is increasing and that

$$0 \leq 1 - x_n = (1-x) \prod_{k=1}^{n} (1 - \pi \gamma_k x_{k-1})$$

$$(21) \qquad \leq (1-x) \exp\left(-\pi \sum_{1\leq k \leq n} \gamma_k x_{k-1}\right)$$

$$(22) \qquad \leq (1-x) \exp(-\pi x \Gamma_n).$$

Plugging (22) into (21) yields

$$0 \leq 1 - x_n \leq (1-x)\exp\left(-\pi\Gamma_n + \pi(1-x)\sum_{1\leq k\leq n}\gamma_k e^{-\pi x \Gamma_{k-1}}\right)$$

$$\leq (1-x)\exp\left(-\pi\Gamma_n + \pi(1-x)e^{\pi x}\int_0^{\Gamma_n} e^{-\pi x u}\,du\right)$$

$$\leq (1-x)\exp\left(\left(\frac{1}{x}-1\right)e^{\pi x}\right)e^{-\pi\Gamma_n}$$

(23) $\qquad = O(e^{-\pi\Gamma_n}).$

On the other hand, for every $n \geq 1$,

$$1 - x_n \geq (1-x)\prod_{k=1}^n (1-\pi\gamma_k).$$

In particular, if one assumes that $\sum_n \gamma_n^2 < +\infty$, then there are some positive real constants $C(x)$ and $C'(x)$ such that, for every $n \geq 1$,

(24) $\qquad C(x) \leq e^{\pi\Gamma_n}(1-x_n) \leq C'(x).$

Now, let us come back to the original procedure with the specification given by (8). By an obvious symmetry argument, one shows, as in the proof of Proposition 4, that the events

$$I_{\infty,x} := \left\{ U_n < 1 - (1-x)\prod_{k=1}^{n-1}(1-\gamma_k \mathbf{1}_{A_k}) \text{ for every } n \geq 1 \right\}$$

and

$$\mathcal{X}^{\uparrow}_{\infty,x} := \left\{ X_n = 1 - (1-x)\prod_{k=1}^{n}(1-\gamma_k \mathbf{1}_{A_k}) \text{ for every } n \geq 0 \right\}$$

satisfy

$$I_{\infty,x} \subset \mathcal{X}^{\uparrow}_{\infty,x}.$$

Now, still following the proof of Proposition 4, $\mathbb{P}_x(I_{\infty,x}) > 0$ as soon as $\sum_n \prod_{k=1}^n (1-p_A\gamma_k) < +\infty$. Moreover, if $\sum_n \gamma_n^2 < +\infty$, the proof of Lemma 2 shows that

$$\frac{\prod_{k=1}^n (1-\gamma_k \mathbf{1}_{A_k})}{\prod_{k=1}^n (1-\gamma_k)^{p_A}} \xrightarrow{\text{a.s.}} \zeta \in (0,+\infty) \qquad \text{as } n \to \infty.$$

Hence,

$$\frac{\prod_{k=1}^n (1-\gamma_k \mathbf{1}_{A_k})}{\prod_{k=1}^n (1-p_A\gamma_k)} \xrightarrow{\text{a.s.}} \zeta' \in (0,+\infty) \qquad \text{as } n \to \infty$$

so that

(25) $$e^{p_A\Gamma_n}\prod_{k=1}^n(1-\gamma_k\mathbf{1}_{A_k})\overset{\text{a.s.}}{\to}\zeta''\in(0,+\infty)\qquad\text{as }n\to\infty.$$

This leads to the following result concerning the rate of convergence of the algorithm (stated here in the infallible case, but an analogous phenomenon occurs in the fallible case for the equilibrium 0).

PROPOSITION 8. *Assume that $0<p_B<p_A\leq 1$ and that*

(26) $$\gamma_n=O(\Gamma_n e^{-p_B\Gamma_n})\quad\text{and}\quad\sum_{n\geq 1}\prod_{k=1}^n(1-p_A\gamma_k)<+\infty.$$

*Then, the two-armed bandit algorithm is a.s. infallible and, for every $x\in(0,1)$, there exists an event of positive $\mathbb{P}_x$-probability $I_{\infty,x}$ on which $X_n$ is nondecreasing and*

(27) $$e^{p_A\Gamma_n}(1-X_n)\overset{\text{a.s.}}{\to}\xi\in(0,+\infty)\qquad\text{as }n\to+\infty.$$

*Assumption* (26) *is fulfilled, for example, when $\gamma_n=\frac{C}{C+n}$, with $\frac{1}{p_A}<C\leq\frac{1}{p_B}$.*

REMARK. (i) Comparing the rates obtained in (24) and in (27), respectively, shows that, for step sequences satisfying (26), the two-armed bandit algorithm converges toward its "target" equilibrium 1 on an event with positive probability infinitely faster than the corresponding algorithm "in average." More generally, the same phenomenon occurs at least at one of the equilibrium points as soon as

$$\sum_n\gamma_n^2<+\infty\quad\text{and}\quad\sum_{n\geq 1}\prod_{k=1}^n(1-\max(p_A,p_B)\gamma_k)<+\infty.$$

This unusual behavior in the field of stochastic approximation is confirmed by the fact that the assumptions of the standard central limit theorem for recursive stochastic algorithms (at rate $\sqrt{\gamma_n}$, see [5] among others) are never fulfilled by the two-armed bandit algorithm: when $p_A\neq p_B$ the martingale increment $\Delta M_n$ involved in the canonical decomposition (6) of the algorithm satisfies

$$\mathbb{E}_x((\Delta M_{n+1})^2/\mathcal{F}_n)\leq X_n(1-X_n)\overset{\text{a.s.}}{\to}0\qquad\text{as }n\to+\infty,$$

whereas this term is supposed to converge toward some positive real number to apply the CLT.

(ii) Proposition 4 can be slightly improved using the same ingredients as above. Namely, if $\sum_n\gamma_n^2<+\infty$, then, for every $x\in(0,1)$,

$$\mathbb{P}_x(\{X_n\text{ goes to }0\text{ monotonously for large enough }n\})>0$$

$$\text{iff }\sum_{n\geq 1}\prod_{k=1}^n(1-p_B\gamma_k)<+\infty$$

and

$$\mathbb{P}_x(\{X_n \text{ goes to } 1 \text{ monotonously for large enough } n\}) > 0$$
$$\text{iff } \sum_{n \geq 1} \prod_{k=1}^{n} (1 - p_A \gamma_k) < +\infty.$$

By symmetry, it suffices to establish the equivalence, for example, for the equilibrium 0. By the Markov property this amounts to showing that if $\sum_n \gamma_n^2 < +\infty$,

$$\mathbb{P}_x(I_{\infty,x}) > 0 \quad \text{if and only if} \quad \sum_{n \geq 1} \prod_{k=1}^{n} (1 - p_A \gamma_k) < +\infty.$$

The equivalence follows from (25) and the Lebesgue dominated convergence Theorem applied to the identity

$$\mathbb{P}_x(I_{\infty,x}) = \mathbb{E}_x\left(\prod_{n \geq 1}\left(1 - \prod_{k=1}^{n}(1 - \gamma_k \mathbf{1}_{A_k})\right)\right).$$

5.2. *Stopping rules.* The proposition below derives an upper-bound for the conditional error probability at time $n$ based on some inequality used in the proof of Lemma 1.

PROPOSITION 9. *Assume that $p_A, p_B \in [0,1]$ and $p_A \neq p_B$. Let $X_\infty = $ a.s.-$\lim_n X_n$ and let $x_\infty = \mathbf{1}_{\{p_A > p_B\}}$ be the "target" parameter of the algorithm. Then, for every $n \geq 1$,*

$$\mathbb{P}_x(X_\infty \neq x_\infty / \mathcal{F}_n)$$
$$\leq \max\left(\min\left(\frac{1 - X_n}{X_n}, \frac{\sum_{k \geq n} \gamma_{k+1}^2}{X_n}\right), \min\left(\frac{X_n}{1 - X_n}, \frac{\sum_{k \geq n} \gamma_{k+1}^2}{1 - X_n}\right)\right).$$

PROOF. Assume for the sake of simplicity that $p_A > p_B$ so that $x_\infty = 1$. Assume that the events $A_n$ and $B_n$ involved in the dynamics of $(X_n)_{n \geq 0}$ are specified by (8). Then, for every $n \geq 1$, one considers $(\bar{X}_k^{(n)})_{k \geq n}$ the (martingale) algorithm defined for every $k \geq n$ by

(28)
$$\bar{X}_n^{(n)} = X_n,$$
$$\bar{X}_{k+1}^{(n)} = \bar{X}_k^{(n)} + \gamma_{k+1} \mathbf{1}_{B_{k+1}}(\mathbf{1}_{\{U_{k+1} \leq \bar{X}_k^{(n)}\}} - \bar{X}_k^{(n)}).$$

It follows from Proposition 3 that, for every $n \geq 0$ and for every $k \geq n$, $\bar{X}_k^{(n)} \leq X_k$, so that $\bar{X}_\infty^{(n)} := \text{a.s.-}\lim_k \bar{X}_k^{(n)} \leq X_\infty$.

Now, as in the proof of Lemma 1, one notices that

$$\mathbb{P}_x(X_\infty = 0/\mathcal{F}_n) \leq \mathbb{P}_x(\bar{X}_\infty^{(n)} = 0/\mathcal{F}_n)$$
$$\leq \frac{\mathbb{E}_x((\bar{X}_\infty^{(n)} - \bar{X}_n^{(n)})^2/\mathcal{F}_n)}{X_n^2}.$$

A straightforward computation based on (28) then shows that the conditional variance increment process of $\bar{X}^{(n)}$ is given for every $k \geq n$ by

$$\langle \bar{X}^{(n)} \rangle_k = p_B \sum_{\ell=n}^{k-1} \gamma_{\ell+1}^2 \bar{X}_\ell^{(n)}(1 - \bar{X}_\ell^{(n)}).$$

Consequently, still as in the proof of Lemma 1,

$$(29) \qquad \mathbb{P}_x(X_\infty = 0/\mathcal{F}_n) \leq \frac{p_B \sum_{k \geq n} \gamma_{k+1}^2 \mathbb{E}_x(\bar{X}_k^{(n)}(1 - \bar{X}_k^{(n)})/\mathcal{F}_n)}{X_n^2}$$
$$\leq \frac{p_B \sum_{k \geq n} \gamma_{k+1}^2 \mathbb{E}_x(\bar{X}_k^{(n)}/\mathcal{F}_n)}{X_n^2}$$
$$= \frac{p_B \bar{X}_n^{(n)} \sum_{k \geq n} \gamma_{k+1}^2}{X_n^2}$$
$$\leq \frac{\sum_{k \geq n} \gamma_{k+1}^2}{X_n}.$$

On the other hand, we know from (7) that

$$\mathbb{E}_x(\bar{X}_k^{(n)}(1 - \bar{X}_k^{(n)})/\mathcal{F}_n) = X_n(1 - X_n) \prod_{\ell=n+1}^{k} (1 - p_B \gamma_\ell^2)$$

so that

$$p_B \sum_{k \geq n} \gamma_{k+1}^2 \mathbb{E}_x(\bar{X}_k^{(n)}(1 - \bar{X}_k^{(n)})/\mathcal{F}_n) = X_n(1 - X_n)\left(1 - \prod_{k \geq n+1} (1 - p_B \gamma_k^2)\right).$$

Plugging this identity in (29) yields

$$\mathbb{P}_x(X_\infty = 0/\mathcal{F}_n) \leq \frac{1 - X_n}{X_n}\left(1 - \prod_{k \geq n+1} (1 - p_B \gamma_k^2)\right) \leq \frac{1 - X_n}{X_n}.$$

The upper-bound for $\mathbb{P}_x(X_\infty = 1/\mathcal{F}_n)$ follows from a symmetry argument.
□

**6. Additional results.**

*Regularity of* $x \mapsto \mathbb{P}_x(X_\infty = 1)$ *when* $p_A > p_B$. One can obtain some regularity results for the function $x \mapsto \mathbb{P}_x(X_\infty = 1)$ as soon as $p_A \neq p_B$ (keep in mind that in that setting, $X_\infty$ is $\{0,1\}$-valued). Namely,

PROPOSITION 10. *If* $p_A > p_B$, *the function* $x \mapsto \mathbb{P}_x(X_\infty = 1)$ *is nondecreasing and analytic on* $(0,1]$.

PROOF. The only point to establish is analyticity. We sketch the proof in the case of a constant step sequence. One starts from the second equality in (11) and the tools developed in the proof of Theorem 2. We also adopt the same notations. Indeed, function $\chi_\gamma$ is analytic on $(0,1]$ since it is an absolutely decreasing function. Then, $\psi_\gamma$ is analytic as well and consequently so is $x \mapsto \mathbb{P}_x(X_n \to 1)$. The extension to nonconstant step sequences is straightforward. □

*About the distribution of* $X_\infty$. When $0 < p_A = p_B \leq 1$ and $\sum_{n \geq 1} \gamma_n^2 < +\infty$, the conditional distribution of $X_\infty$ given $\{X_\infty \neq 0, 1\}$ is continuous. This follows from Theorem 3.IV.13 in [5].

*Still open questions...* The main open question is, of course, to find a necessary and sufficient condition for the algorithm to be a.s. infallible. For example, when $p_A = p_B = 1$, assumption (3) is easier to express using the partial sums $S_n$ of the $\Delta_n$'s by

$$\tag{30} \sum_{n \geq 1} \frac{1}{S_n} < +\infty.$$

If $\Delta_n = \log n \log_2^\beta n$, assumption (30) is equivalent to $\beta > 1$, whereas $\Delta_n = O(\Gamma_n)$ in Lemma 1 reads $\beta = 0$. So we are facing a log log problem.

Furthermore, it follows from the Borel–Cantelli lemma for independent events that

$$\limsup_n \frac{Y_n}{\Delta_n} \geq \limsup_n \mathbf{1}_{\{U_n \geq x/S_{n-1}\}} = 1, \qquad \mathbb{P}_x\text{-a.s.}$$

when $\sum_{n \geq 1} 1/S_n = +\infty$.

It is to be noticed that, when $\Delta_n = \log n \log_2^\beta n$ for some $\beta \in (0,1)$, this straightforwardly implies that $\limsup_n Y_n / \log S_n = +\infty$ (which was the key step of Lemma 1). Unfortunately, for such sequences $\Delta_n$, Lemma 1 only implies that $\mathbb{P}_x(X_\infty = 0) = 0$ iff $\limsup_n \frac{Y_n}{\Delta_n} > 1$, $\mathbb{P}_x$-a.s.

*The last remark.* Let $\sigma_n := \min\{k > \sigma_{n-1} / U_k < X_{k-1}\}$ and $\sigma_0 := 0$ denote the increasing break times. Assumption (30) is equivalent to $\mathbb{P}_x(\sigma_1 = +\infty) > 0$ for every $x \in (0,1)$.

Otherwise, all the $\sigma_n$'s are $\mathbb{P}_x$-a.s. finite for every $x \in (0,1)$: this follows from the expression of $\mathbb{P}_x(\sigma_1 \geq k)$ and from the Markov property.

**Acknowledgment.** The authors are grateful to Jean-Claude Fort for stimulating discussions.

D. Lamberton
Laboratoire d'analyse
  et de mathématiques appliquées
UMR 8050
Université Marne-la-Vallée
Cité Descartes 5, Bld Descartes
Champs-sur-Marne
F-77454 Marne-la-Vallée Cedex 2
France
e-mail: dlamb@math.univ-mlv.fr

G. Pagès
Laboratoire de Probabilités
  et Modélisation aléatoire
UMR 7599
Université Paris 6
Case 188
4 pl. Jussieu
F-75252 Paris Cedex 5
France
e-mail: gpa@ccr.jussieu.fr

P. Tarrès
Laboratoire de Statistique
  et Probabilités
UMR CNRS C5583
Université Paul Sabatier
118 route de Narbonne
F-31062 Toulouse Cedex 4
France
e-mail: tarres@math.ups-tlse.fr